\documentclass[leqno,12pt]{article}
\usepackage{amsmath, amssymb}
\usepackage{amsthm}

\theoremstyle{plain}
\newtheorem{theorem}{\indent\sc Theorem}[section]
\newtheorem{lemma}[theorem]{\indent\sc Lemma}
\newtheorem{corollary}[theorem]{\indent\sc Corollary}

\theoremstyle{definition}
\newtheorem{Def}[theorem]{\indent\sc Definition}
\newtheorem{remark}[theorem]{\indent\sc Remark}



\newcommand{\field}[1]{\mathbb{#1}}

\newcommand{\R}{\field{R}}

\newcommand{\beeq}{\begin{equation}}
\newcommand{\eneq}{\end{equation}}

\newcommand{\ba}{\begin{array}}
\newcommand{\ea}{\end{array}}

\newcommand{\be}{\begin{equation}}
\newcommand{\ee}{\end{equation}}
\newcommand{\bea}{\begin{eqnarray}}
\newcommand{\eea}{\end{eqnarray}}

\begin{document}

{\centerline{ \bf{On the regularity of weak solutions of the Boussinesq equations in
Besov spaces}}}
$$\,$$
{\centerline{ \bf{Dedicated to Enrique Zuazua on the occasion of his sixtieth birthday}}}
$$\,$$
{\center{Annamaria Barbagallo\\
Department of Mathematics and Applications "R.Caccioppoli"\\University of Naples "Federico II", Naples
e-mail:annamaria.barbagallo@unina.it
}
\center{
Sadek Gala\\
Department of Sciences exactes, Ecole Normale Sup´erieure de Mostaganem
Box 227, Mostaganem 27000, Algeria
e-mail: sgala793@gmail.com
}
\center{
Maria Alessandra Ragusa\\
Dipartimento di Matematica e Informatica, Universit\`{a} di Catania,
Viale Andrea Doria, 6-95125 Catania, Italy,\\
e-mail:maragusa@dmi.unict.it
}
\center{
 Michel Th\`{e}ra\\
XLIM UMR-CNRS 7252 Universit\`{e} de Limoges and Centre for Informatics and Applied Optimisation,
Federation University Australia
e-mail: michel.thera@unilim.fr
}}

$$\,$$


{\centerline{\bf{ Abstract}}}

The main issue addressed in this paper concerns  an extension of a result by Z. Zhang
who proved,   in the context of the homogeneous Besov space $\dot{B}_{\infty ,\infty }^{-1}(\mathbb{R}%
^{3})$,   that, if the solution of the Boussinesq equation (\ref%
{eq1.1}) below (starting with an initial data in $H^{2}$) is such that $%
(\nabla u,\nabla \theta )\in L^{2}\left( 0,T;\dot{B}_{\infty ,\infty }^{-1}(%
\mathbb{R}^{3})\right)$,
then the solution remains smooth
forever after $T$. In this contribution, we prove the same result for weak solutions just by
assuming the condition on the velocity $u$ and not on the temperature $\theta$.
$$\,$$
\noindent {\textit{Mathematics Subject Classification(2000):}\thinspace
\thinspace 35Q35,\thinspace \thinspace 35B65,\thinspace \thinspace 76D05.}


\section{Introduction}

We are interested in the regularity of weak solutions of the Cauchy problem
related to the Boussinesq equations in $\mathbb{R}^{3}:$
\begin{equation}
\left\{
\begin{array}{c}
\partial _{t}u+\left( u\cdot \nabla \right) u-\Delta u+\nabla \pi =\theta
e_{3}, \\
\partial _{t}\theta +\left( u\cdot \nabla \right) \theta -\Delta \theta =0,
\\
\nabla \cdot u=0, \\
u(x,0)=u_{0}(x),\text{ \ \ }\theta (x,0)=\theta _{0}(x),%
\end{array}%
\right.  \label{eq1.1}
\end{equation}%
where $x\in \mathbb{R}^{3}$ and $t\geq 0$.\   Here,  $u:\R^3\times \R_+\to \R^3$  
 is the velocity
field of the flow, $\pi =\pi (x,t)\in \mathbb{R}$ is a scalar function
representing the pressure,  $\theta:\R^3\times \R_+\to \R^3$
represents the
temperature of the fluid  and $e_{3}=(0,0,1)^{T}$. Note that  $u_{0}(x)$ and $\theta _{0}(x)$
are given initial velocity and initial temperature with $\nabla \cdot
u_{0}=0 $ in the sense of distributions.

Owing to the physical importance and the mathematical challenges, the study
of (\ref{eq1.1})   {which  describes the dynamics of a viscous incompressible fluid with heat exchanges,
has  a long history and has attracted many contributions from physicists and mathematicians \cite{Maj}.
Although  Boussinesq  equations consist  in a simplification of the original 3-D incompressible flow, they share a similar vortex stretching effect. For this reason they retains most of  the mathematical and physical difficulties of the  3-D incompressible flow, and  therefore,  these equations have been studied and applied to various  fields. Examples include    for instance   geophysical applications, where  they  serve as a model, see, e.g. \cite{Pedlosky}. There are several other results on existence and blowup criteria in different kinds of spaces which have been obtained, (see \cite{Abidi,Cannon, Danchin08,Danchin, Ye2}).

The problem of the global-in-time well-posedness of (\ref{eq1.1}) in a
three-dimensional space is highly challenging, due to the fact that the
system contains the incompressible 3D Navier-Stokes equations as a special case (obtained
by setting $\theta =0$), for which the issue of global well-posedness has
not been proved until now. However, the question of the regularity of weak
solutions is an outstanding open problem in mathematical fluid mechanics and
many interesting results have been obtained (see e.g. \cite{CN, CKN, FO, FZ,
GGR, GG, GR, QYW, QDY, X, XZZ, Ye2, Ye}). We are interested in the classical
problem of finding sufficient conditions for weak solutions of (\ref{eq1.1})
such that they become regular.

Realizing the dominant role played by the velocity field in the regularity
issue, Ishimura and Morimoto  \cite{IM} were able to derive criteria in terms of
the velocity field $u$ alone. They showed that, if $u$ satisfies%
\begin{equation}
\nabla u\in L^{1}\left( 0,T;L^{\infty }(\mathbb{R}^{3})\right) ,
\label{eq210}
\end{equation}%
then the solution $(u,\theta )$ is regular on $[0,T]$. It is worthy to
emphasize that there are no assumptions on the temperature $\theta $. This
assumption (\ref{eq210}) was weakened in [6] with the $L^{\infty }-$norm
replaced by norms in Besov spaces $\dot{B}_{\infty ,\infty }^{0}$. Quite
recently, Z.\ Zhang  \cite{Z} showed that $(u,\theta )$ is a strong solution if
\begin{equation}
(\nabla u,\nabla \theta )\in L^{2}\left( 0,T;\dot{B}_{\infty ,\infty }^{-1}(%
\mathbb{R}^{3})\right) ,  \label{eq2}
\end{equation}%
where $\dot{B}_{\infty ,\infty }^{-1}$ denotes the homogenous Besov space. A
logarithmically improvement of Zhang's result, controlled by its $H^{3}-$%
norm, was given by Ye \cite{Ye2}.

The main purpose of this work is to establish an improvement of Zhang's
regularity criterion (\ref{eq2}). Now,  the refined regularity criterion in
terms of the gradient of the velocity $\nabla u$ can be stated as follows:

\begin{theorem}[Main result] 
\label{th1}
Assume that $\left( u_{0},\theta _{0}\right) \in L^{2}(\mathbb{R}%
^{3})$ with $\nabla \cdot u_{0}=0$. Let $(u,\theta )$ be a weak solution to the
Boussinesq equations on some interval $(0,T)$ with $0<T\leq \infty $.\ If
\begin{equation}
\nabla u\in L^{2}\left( 0,T;\dot{B}_{\infty ,\infty }^{-1}(\mathbb{R}%
^{3})\right) ,  \label{eq2.5}
\end{equation}%
then the weak solution $(u,\theta )$ is regular in $(0,T]$, that is $%
(u,\theta )\in C^{\infty }\left( \mathbb{R}^{3}\times (0,T]\right) $.
\end{theorem}

\begin{remark}
This result is expected because of the fact that the (refinement of)
Beale-Kato-Majda type criterion is well known in the class $\dot{B}_{\infty
,\infty }^{0}$ for the 3D Boussinesq equations and one may replace the
vorticity by $\nabla u$ since the Riesz transforms are continuous in $\dot{B}%
_{\infty ,\infty }^{0}$. Then, the temperature plays a less dominant role
than the velocity field does in the regularity theory of solutions to
the Boussinesq equations. Furthermore, clearly  Theorem \ref{th1}  is
an improvement of Zhang's regularity criterion (\ref{eq2}).
\end{remark}

By a weak solution, we mean that $(u,\theta ,\pi )$ must satisfy (\ref{eq1.1}%
) in the sense of  distributions. In addition, we have the basic regularity
for the weak solution%
\begin{equation*}
(u,\theta )\in L^{\infty }(0,T;L^{2}(\mathbb{R}^{3}))\cap L^{2}(0,T;H^{1}(%
\mathbb{R}^{3})),
\end{equation*}%
for any $T>0$. If a weak solution $(u,\theta )$ satisfies
\begin{equation*}
(u,\theta )\in L^{\infty }(0,T;H^{1}(\mathbb{R}^{3}))\cap L^{2}(0,T;H^{2}(%
\mathbb{R}^{3})),
\end{equation*}%
then actually $(u,\theta )$ is a strong (classical) solution. It is worth to
note that for strong solutions, we can gain more regularity properties.

Throughout this paper, $C$ denotes a generic positive constant which may
 vary from one line to another.

\section{Preliminaries}

In this section we introduce the function spaces that will be used to
state and prove the main result, and we collect and/or derive a number of
auxiliary estimates that will be needed throughout the proof. Before
introducing the homogeneous Besov and Triebel-Lizorkin spaces, we have to
fix some notations. By $\mathcal{S}$ we denote the class of rapidly
decreasing functions. The dual space of $\mathcal{S}$, i.e., the space of
tempered distributions on $\mathbb{R}^{3}$ is denoted by $\mathcal{S}%
^{\prime }$. For $u\in \mathcal{S}(\mathbb{R}^{3})$, the Fourier transform
of $u$ is defined by%
\begin{equation*}
\mathcal{F}u(\omega )=\widehat{u}(\omega )=\int_{\mathbb{R}%
^{3}}u(x)e^{-ix\cdot \omega }dx,\text{ \ \ }\omega \in \mathbb{R}^{3}.
\end{equation*}%
The homogeneous Littlewood-Paley decomposition relies upon a dyadic
partition of unity. We can use for instance any $\varphi \in \mathcal{S}(%
\mathbb{R}^{3})$, supported in $\mathcal{C}\triangleq \left\{ \omega \in
\mathbb{R}^{3}:\frac{3}{4}\leq \left\vert \omega \right\vert \leq \frac{8}{3}%
\right\} $ such that
\begin{equation*}
\sum\limits_{l\in \mathbb{Z}}\varphi (2^{-l}\omega )=1\text{ \ \ if \ }%
\omega \neq 0.\text{ }
\end{equation*}%
Denoting $h=\mathcal{F}^{-1}\varphi $, we then define dyadic blocks  in this way:%
\begin{equation*}
\Delta _{l}u\triangleq \varphi (2^{-l}D)u=2^{3l}\int_{\mathbb{R}%
^{3}}h(2^{l}y)u(x-y)dy,\text{ \ \ for each \ }l\in \mathbb{Z},
\end{equation*}%
and%
\begin{equation*}
S_{l}u\triangleq \sum\limits_{k\leq l-1}\Delta _{k}u.
\end{equation*}%
The formal decomposition%
\begin{equation*}
u=\sum\limits_{l\in \mathbb{Z}}\Delta _{l}u
\end{equation*}%
is called the homogeneous Littlewood-Paley decomposition.

\begin{remark}
The above dyadic decomposition has nice properties of quasi-orthogonality:
with  our choice of $\varphi $, we have,
\begin{equation*}
\Delta _{k}\Delta _{l}u\equiv 0\text{ \ if \ }\left\vert k-l\right\vert \geq
2\text{ \ \ \ and \ \ }\Delta _{k}(S_{k-1}u\Delta _{l}u)\equiv 0\text{ \ if
\ }\left\vert k-l\right\vert \geq 5.
\end{equation*}
\end{remark}

With the introduction of $\Delta _{l}$, let us recall the definition of
homogeneous Besov and Triebel-Lizorkin spaces (see \cite{Tri} for more
details).

\begin{Def}
The homogeneous Besov space $\dot{B}_{p,q}^{s}(\mathbb{R}^{3})$ is defined
by
\begin{equation*}
\dot{B}_{p,q}^{s}(\mathbb{R}^{3})=\left\{ u\in \mathcal{S}^{\prime }(\mathbb{%
R}^{3})/\mathcal{P}(\mathbb{R}^{3}):\left\Vert u\right\Vert _{\dot{B}%
_{p,q}^{s}}<\infty \right\} ,
\end{equation*}%
for $s\in \mathbb{R}$ and $1\leq p,q\leq \infty $, where
\begin{equation*}
\left\Vert u\right\Vert _{\dot{B}_{p,q}^{s}}=\left\{
\begin{array}{c}
\left( \sum\limits_{j\in \mathbb{Z}}2^{jsq}\left\Vert \Delta
_{j}u\right\Vert _{L^{p}}^{q}\right) ^{\frac{1}{q}},\text{ \ \ if \ \ \ }%
1\leq q<\infty , \\
\underset{j\in \mathbb{Z}}{\sup }2^{js}\left\Vert \Delta _{j}u\right\Vert
_{L^{p}},\text{ \ \ if \ \ \ }q=\infty ,%
\end{array}%
\right.
\end{equation*}%
and $\mathcal{P}(\mathbb{R}^{3})$  is  the set of all scalar polynomials
defined on $\mathbb{R}^{3}$. Similarly, the homogeneous Triebel-Lizorkin
spaces $\dot{F}^s_{p,q}(\mathbb{R}^{3})$ is a quasi-normed
space equipped with the family of semi-norms $\left\Vert \cdot \right\Vert _{%
\dot{F}^s_{p,q}}$ which are defined by
\begin{equation*}
\left\Vert u\right\Vert_{\dot{F}^s_{p,q}}=\left\{
\begin{array}{c}
\left\Vert \left( \sum\limits_{j\in \mathbb{Z}}2^{jsq}\left\vert \Delta
_{j}u\right\vert ^{q}\right) ^{\frac{1}{q}}\right\Vert _{L^{p}},\text{ \ \ \
if \ \ }1\leq q<\infty , \\
\left\Vert \underset{j\in \mathbb{Z}}{\sup }2^{js}\left\vert \Delta
_{j}u\right\vert \right\Vert _{L^{p}},\text{ \ \ \ if \ \ }q=\infty .%
\end{array}%
\right.
\end{equation*}
\end{Def}

Notice that there exists a universal constant $C$ such that
\begin{equation*}
C^{-1}\left\Vert u\right\Vert _{\dot{B}_{p,q}^{s}}\leq \left\Vert \nabla
u\right\Vert _{\dot{B}_{p,q}^{s-1}}\leq C\left\Vert u\right\Vert _{\dot{B}%
_{p,q}^{s}}.
\end{equation*}%
In particular,
\begin{equation*}
u\in \dot{B}_{\infty ,\infty }^{0}\left( \mathbb{R}^{3}\right)
\Longleftrightarrow \nabla u\in \dot{B}_{\infty ,\infty }^{-1}\left( \mathbb{%
R}^{3}\right) .
\end{equation*}%
From this observation we derive the following corollary to Theorem \ref{th1}.

\begin{corollary}
Suppose that $\left( u,\theta \right) $ is a weak solution of the  Boussinesq
equations on $(0,T)$.\ If
\begin{equation}
u\in L^{2}\left( 0,T;\dot{B}_{\infty ,\infty }^{0}\left( \mathbb{R}%
^{3}\right) \right) ,  \label{eq2.6}
\end{equation}%
then the weak solution $\left( u,\theta \right) $ is regular in $(0,T].$
\end{
corollary}

Next, we introduce the following Bernstein lemma due to \cite{Che}.

\begin{lemma}[Bernstein]
For all $k\in \mathbb{N}$, $j\in \mathbb{Z}$, and $1\leq p\leq q\leq \infty$, we have for all $f\in \mathcal{S}(\mathbb{R}^{3})$ :

\begin{description}
\item[(i)]
\begin{equation*}
\underset{\left\vert \alpha \right\vert =k}{\sup }\left\Vert \nabla ^{\alpha
}\Delta _{j}f\right\Vert _{L^{q}}\leq C_{1}2^{jk+3j(\frac{1}{p}-\frac{1}{q}%
)}\left\Vert \Delta _{j}f\right\Vert _{L^{p}}
\end{equation*}

\item[(ii)]
\begin{equation*}
\left\Vert \Delta _{j}f\right\Vert _{L^{p}}\leq C_{2}2^{-jk}\underset{%
\left\vert \alpha \right\vert =k}{\sup }\left\Vert \nabla ^{\alpha }\Delta
_{j}f\right\Vert _{L^{p}},
\end{equation*}%
where $C_{1}$, $C_{2}$ are positive constants independent of $f$ and $j$.
\end{description}
\end{
lemma}

The proof of the main result needs a logarithmic Sobolev inequality in terms
of Besov space. It will  play  an important role in the proof of Theorem \ref%
{th1}. The following is a well-known embedding result, (cf. \cite{Tri}, pp.
244):
\begin{equation*}
L^{\infty }(\mathbb{R}^{3})\hookrightarrow BMO(\mathbb{R}^{3})=\overset{%
\cdot }{F}_{\infty ,2}^{0}\hookrightarrow \dot{B}_{\infty ,\infty }^{0}(%
\mathbb{R}^{3}),
\end{equation*}%
where $BMO(\mathbb{R}^{3})$ stands for  the Bounded Mean Oscillations space \cite%
{Tri}.
\vskip 2mm
We state and prove the following lemma.

\begin{
lemma}

\label{lem1}Suppose that $\nabla f\in \dot{B}_{\infty ,\infty }^{-1}(\mathbb{%
R}^{3})$ and $f\in H^{s}(\mathbb{R}^{3})$ for all $s>\frac{3}{2}$. Then,
there exists a constant $C>0$ such that
\begin{equation}
\left\Vert f\right\Vert _{L^{\infty }}\leq C\left[ 1+\left\Vert \nabla
f\right\Vert _{\dot{B}_{\infty ,\infty }^{-1}}\left( \ln ^{+}\left\Vert
f\right\Vert _{H^{s}}\right) ^{\frac{1}{2}}\right] ,  \label{eq21}
\end{equation}%
holds, where $H^{s}$\ denotes the standard Sobolev space and
\begin{equation*}
\ln ^{+}x=\left\{
\begin{array}{c}
\ln x,\text{ \ \ if \ \ }x>e, \\
1,\text{ \ \ if \ \ }0<x\leq e.%
\end{array}%
\right.
\end{equation*}
\end{
lemma}
\begin{
proof}

The proof is an easy modification of the one in \cite{KOT}. Owing the
Littlewood--Paley decomposition, we can rewrite%
\begin{equation*}
f=\sum\limits_{j\in \mathbb{Z}}\Delta _{j}f=\sum\limits_{j=-\infty
}^{-N-1}\Delta _{j}f+\sum\limits_{j=-N}^{N}\Delta
_{j}f+\sum\limits_{j=N+1}^{+\infty }\Delta _{j}f,
\end{equation*}%
where $N$ is a positive integer to be determined later. Bernstein's lemma
and Young's inequality give rise to%
\begin{eqnarray}
\left\Vert f\right\Vert _{L^{\infty }} &\leq &\sum\limits_{j=-\infty
}^{-N-1}\left\Vert \Delta _{j}f\right\Vert _{L^{\infty
}}+\sum\limits_{j=-N}^{N}\left\Vert \Delta _{j}f\right\Vert _{L^{\infty
}}+\sum\limits_{j=N+1}^{+\infty }\left\Vert \Delta _{j}f\right\Vert
_{L^{\infty }}  \notag \\
&\leq &C\sum\limits_{j<-N}2^{\frac{3}{2}j}\left\Vert \Delta _{j}f\right\Vert
_{L^{2}}+CN\left\Vert f\right\Vert _{\dot{B}_{\infty ,\infty
}^{0}}+C\sum\limits_{j>N}2^{(-s+\frac{3}{2})j}\left\Vert \Delta
_{j}f\right\Vert _{L^{2}}2^{js}  \notag \\
&\leq &C\left( 2^{-\frac{3}{2}N}\left\Vert f\right\Vert _{L^{2}}+N\left\Vert
\nabla f\right\Vert _{\dot{B}_{\infty ,\infty
}^{-1}}+\sum\limits_{j>N}2^{(-s+\frac{3}{2})j}\left\Vert f\right\Vert _{\dot{%
B}_{2,\infty }^{s}}\right)   \notag \\
&\leq &C\left( 2^{-\frac{3}{2}N}\left\Vert f\right\Vert _{L^{2}}+N\left\Vert
\nabla f\right\Vert _{\dot{B}_{\infty ,\infty }^{-1}}+2^{(-s+\frac{3}{2}%
)N}\left\Vert f\right\Vert _{H^{s}}\right) ,  \label{eq19}
\end{eqnarray}%
where we have used the fact that $s>\frac{3}{2}$ and the Besov embedding $%
H^{s}\hookrightarrow \dot{B}_{2,\infty }^{s}$.

Setting $\alpha =\min \left( s-\frac{3}{2},\frac{3}{2}\right) $, we derive
\begin{equation}\label{taormina}
\left\Vert f\right\Vert _{L^{\infty }}\leq C\left( 2^{-\alpha N}\left\Vert
f\right\Vert _{H^{s}}+N\left\Vert \nabla f\right\Vert _{\dot{B}_{\infty
,\infty }^{-1}}\right).
\end{equation}

Now choose $N$ such that $2^{-\alpha N}\left\Vert f\right\Vert
_{H^{s}}\leq 1$. Thus we get%

$N\geq \frac{\log \left\Vert f\right\Vert _{H^{s}}}{\alpha \log 2}.$
\end{
proof}
Next, the following lemma is needed.

\begin{
lemma}
\label{lem2}Let $g,h\in H^{1}(\mathbb{R}^{3})$ and $f\in BMO(\mathbb{R}^{3})$%
. Then we have%
\begin{equation}
\int_{\mathbb{R}^{3}}f\cdot \nabla (gh)dx\leq C\left\Vert f\right\Vert
_{BMO}\left( \left\Vert \nabla g\right\Vert _{L^{2}}\left\Vert h\right\Vert
_{L^{2}}+\left\Vert g\right\Vert _{L^{2}}\left\Vert \nabla h\right\Vert
_{L^{2}}\right) .  \label{eq17}
\end{equation}
\end{
lemma}
\begin{
proof}

The proof of the above lemma requires some paradifferential calculus. We
have to recall here that paradifferential calculus enables to define a
generalized product between distributions. It is continuous in many
functional spaces where the usual product does not make sense (see the
pioneering work of J.-M. Bony in \cite{Bony}). The paraproduct between $f$
and $g$ is defined by%
\begin{equation*}
T_{f}g\triangleq \sum\limits_{j\in \mathbb{Z}}S_{j-1}f\Delta _{j}g.
\end{equation*}%
We thus have the following formal decomposition (modulo a polynomial):
\begin{equation}
fg=T_{f}g+T_{g}f+R(f,g),  \label{eq18}
\end{equation}%
with%
\begin{equation*}
R(f,g)=\sum\limits_{\left\vert j-k\right\vert \leq 1}\Delta _{j}f\Delta
_{k}g.
\end{equation*}%
Coming back to the proof of Lemma \ref{lem2}, we split $\int_{\mathbb{R}%
^{3}}f\cdot \nabla (gh)dx$ into%
\begin{eqnarray*}
\int_{\mathbb{R}^{3}}f\cdot \nabla (gh)dx &=&\int_{\mathbb{R}^{3}}f\cdot
\nabla (T_{g}h)dx+\int_{\mathbb{R}^{3}}f\cdot \nabla (gT_{h})dx+\int_{%
\mathbb{R}^{3}}f\cdot \nabla R(g,h)dx \\
&=&I_{1}+I_{2}+I_{3}.
\end{eqnarray*}%
Since we know that $BMO=\overset{\cdot }{F}_{\infty ,2}^{0}$ (see pp.
243--244 of \cite{Tri}), the duality  between $\overset{\cdot }{F}_{\infty
,2}^{0}$ and $\overset{\cdot }{F}_{1,2}^{0}$ guarantees that
\begin{eqnarray*}
I_{1} &=&\int_{\mathbb{R}^{3}}f\cdot (T_{\nabla g}h)dx+\int_{\mathbb{R}%
^{3}}f\cdot (T_{g}\nabla h)dx \\
&\leq &\left\Vert f\right\Vert _{BMO}(\left\Vert T_{\nabla g}h\right\Vert _{%
\overset{\cdot }{F}_{1,2}^{0}}+\left\Vert T_{g}\nabla h\right\Vert _{\overset%
{\cdot }{F}_{1,2}^{0}}) \\
&=&\left\Vert f\right\Vert _{BMO}(I_{11}+I_{12}).
\end{eqnarray*}%
In view of the boundedness of the Hardy-Littlewood maximal operator $%
\mathcal{M}$ in $L^{p}$ spaces $(1<p<\infty )$ (c.f. Stein [\cite{St}, Chap
II, Theorem 1]), we can estimate the term $I_{11}$ as follows :%
\begin{eqnarray*}
I_{11} &\approx &\left\Vert \left( \sum\limits_{j\in \mathbb{Z}}\left\vert
S_{j-1}(\nabla g)\right\vert ^{2}\left\vert \Delta _{j}h\right\vert
^{2}\right) ^{\frac{1}{2}}\right\Vert _{L^{1}}\leq C\left\Vert \mathcal{M}%
(\nabla g)\left( \sum\limits_{j\in \mathbb{Z}}\left\vert \Delta
_{j}h\right\vert ^{2}\right) ^{\frac{1}{2}}\right\Vert _{L^{1}} \\
&\leq &C\left\Vert \mathcal{M}(\nabla g)\right\Vert _{L^{2}}\left\Vert
\left( \sum\limits_{j\in \mathbb{Z}}\left\vert \Delta _{j}h\right\vert
^{2}\right) ^{\frac{1}{2}}\right\Vert _{L^{2}}\leq C\left\Vert \nabla
g\right\Vert _{L^{2}}\left\Vert h\right\Vert _{L^{2}}.
\end{eqnarray*}%
Repeating the same arguments, we also have for $I_{12}$%
\begin{eqnarray*}
I_{12} &\approx &\left\Vert \left( \sum\limits_{j\in \mathbb{Z}}\left\vert
S_{j-1}(g)\right\vert ^{2}\left\vert \Delta _{j}(\nabla h)\right\vert
^{2}\right) ^{\frac{1}{2}}\right\Vert _{L^{1}}\leq C\left\Vert \mathcal{M}%
(g)\left( \sum\limits_{j\in \mathbb{Z}}\left\vert \Delta _{j}(\nabla
h)\right\vert ^{2}\right) ^{\frac{1}{2}}\right\Vert _{L^{1}} \\
&\leq &C\left\Vert g\right\Vert _{L^{2}}\left\Vert \left( \sum\limits_{j\in
\mathbb{Z}}\left\vert \Delta _{j}(\nabla h)\right\vert ^{2}\right) ^{\frac{1%
}{2}}\right\Vert _{L^{2}}\leq C\left\Vert g\right\Vert _{L^{2}}\left\Vert
\nabla h\right\Vert _{L^{2}}.
\end{eqnarray*}%
Collecting these estimates, we obtain
\begin{equation*}
I_{1}\leq C\left\Vert f\right\Vert _{BMO}\left( \left\Vert \nabla
g\right\Vert _{L^{2}}\left\Vert h\right\Vert _{L^{2}}+\left\Vert
g\right\Vert _{L^{2}}\left\Vert \nabla h\right\Vert _{L^{2}}\right) .
\end{equation*}%
As a result, estimating $I_{2}$ following the same arguments, we obtain
\begin{equation*}
I_{2}\leq C\left\Vert f\right\Vert _{BMO}\left( \left\Vert \nabla
g\right\Vert _{L^{2}}\left\Vert h\right\Vert _{L^{2}}+\left\Vert
g\right\Vert _{L^{2}}\left\Vert \nabla h\right\Vert _{L^{2}}\right) .
\end{equation*}%
For the third term $I_{3}$, using the embedding relation $\dot{B}%
_{1,1}^{0}\subset \overset{\cdot }{F}_{1,2}^{0}$ and in view of Bernstein's
lemma, we can deduce that%
\begin{eqnarray*}
I_{3} &\leq &\left\Vert f\right\Vert _{BMO}\left\Vert \nabla
R(g,h)\right\Vert _{\overset{\cdot }{F}_{1,2}^{0}}\leq C\left\Vert
f\right\Vert _{BMO}\left\Vert R(g,h)\right\Vert _{\dot{B}_{1,1}^{1}} \\
&\leq &C\left\Vert f\right\Vert _{BMO}\sum\limits_{j\in \mathbb{Z}%
}2^{j}\left\Vert \Delta _{j}g\cdot \widetilde{\Delta }_{j}h\right\Vert
_{L^{1}} \\
&\leq &C\left\Vert f\right\Vert _{BMO}\sum\limits_{j\in \mathbb{Z}%
}2^{j}\left\Vert \Delta _{j}g\right\Vert _{L^{2}}\left\Vert \widetilde{%
\Delta }_{j}h\right\Vert _{L^{2}} \\
&\leq &C\left\Vert f\right\Vert _{BMO}\left\Vert \nabla g\right\Vert
_{L^{2}}\left\Vert h\right\Vert _{L^{2}}.
\end{eqnarray*}%
so that the proof of Lemma \ref{lem2} is achieved.
\end{
proof}
We often use the following well-known lemma.

\begin{
lemma}
\label{lem3} Let $1\leq q,r<\infty $ and $m\leq k$. Suppose that $\theta $
and $j$ satisfy $m\leq j\leq k$, $0\leq \theta \leq 1$ and define $p\in
\lbrack 1,+\infty ]$ by
\begin{equation*}
\frac{1}{p}=\frac{j}{3}+\theta (\frac{1}{r}-\frac{m}{3})+(1-\theta )(\frac{1%
}{q}-\frac{k}{3}).
\end{equation*}%
Then, the inequality
\begin{equation*}
\left\Vert \nabla ^{j}f\right\Vert _{L^{p}}\leq C\left\Vert \nabla
^{m}f\right\Vert _{L^{q}}^{1-\theta }\left\Vert \nabla ^{k}f\right\Vert
_{L^{r}}^{\theta }\text{ \ \ \ \ for \ \ }f\in W^{m,q}(\mathbb{R}^{3})\cap
W^{k,r}(\mathbb{R}^{3}),
\end{equation*}%
holds with some constant $C>0$.
\end{
lemma}

\section{Proof of Theorem \protect\ref{th1}}

\bigskip Now we are ready to prove our main result of this section.
\vskip 2mm
{\bf Proof}

First, note that a weak solution $(u,\theta )$ to (\ref{eq1.1}) has at least
one global weak solution%
\begin{equation*}
(u,\theta )\in L^{\infty }(0,T;L^{2}(\mathbb{R}^{3}))\cap L^{2}(0,T;H^{1}(%
\mathbb{R}^{3})),
\end{equation*}%
which satisfies the following energy inequality%
\begin{eqnarray*}
&&\frac{1}{2}(\left\Vert u(\cdot ,t)\right\Vert _{L^{2}}^{2}+\left\Vert
\theta (\cdot ,t)\right\Vert _{L^{2}}^{2})+\int_{0}^{t}(\left\Vert \nabla
u(\cdot ,\tau )\right\Vert _{L^{2}}^{2}+\left\Vert \nabla \theta (\cdot
,\tau )\right\Vert _{L^{2}}^{2})d\tau \\
&\leq &\frac{1}{2}(\left\Vert u_{0}\right\Vert _{L^{2}}^{2}+\left\Vert
\theta _{0}\right\Vert _{L^{2}}^{2}),
\end{eqnarray*}%
for almost every $t\geq 0$.

In order to prove that  $(u,\theta )\in C^{\infty }\left( \mathbb{R}^{3}\times
(0,T]\right) $, as it is well known, it suffices to show that the weak
solution $(u,\theta )$ is also a strong solution on $(0,T]$, which means
that:%
\begin{equation*}
(u,\theta )\in L^{\infty }(0,T;H^{1}(\mathbb{R}^{3}))\cap L^{2}(0,T;H^{2}(%
\mathbb{R}^{3})).
\end{equation*}%
Owing to (\ref{eq2.5}), we know that for any small constant $\epsilon >0$,
there exists $T_{0}=T_{0}(\epsilon )<T$ such that%
\begin{equation}
\int_{T_{0}}^{T}\left\Vert \nabla u(\cdot ,\tau )\right\Vert _{\dot{B}%
_{\infty ,\infty }^{-1}}^{2}d\tau \leq \epsilon .  \label{eq51}
\end{equation}%
To do so, we shall work on the local strong solution with the initial datum $%
\left( u_{0},\theta _{0}\right) $ on its maximal existence time interval $%
(0,T_{0})$. Then, we have only to show that%
\begin{equation*}
\underset{0\leq t<T_{0}}{\sup }(\left\Vert \nabla u(\cdot ,t)\right\Vert
_{L^{2}}^{2}+\left\Vert \nabla \theta (\cdot ,t)\right\Vert
_{L^{2}}^{2})+\int_{0}^{T_{0}}(\left\Vert \nabla u(\cdot ,\tau )\right\Vert
_{L^{2}}^{2}+\left\Vert \nabla \theta (\cdot ,\tau )\right\Vert
_{L^{2}}^{2})d\tau \leq C<\infty ,
\end{equation*}%
here and in what follows $C$ denotes various positive constants which are
independent from $T_{0}$.

Take the operator $\nabla $ in equations (\ref{eq1.1})$_{1}$ and (\ref{eq1.1}%
)$_{2}$, respectively, and the scalar product of them $\nabla u$ and $\nabla
\theta $, respectively and add them together, to  obtain%
\begin{eqnarray}
&&\frac{1}{2}\frac{d}{dt}(\left\Vert \nabla u\right\Vert
_{L^{2}}^{2}+\left\Vert \nabla \theta \right\Vert _{L^{2}}^{2})+\left\Vert
\Delta u\right\Vert _{L^{2}}^{2}+\left\Vert \Delta \theta \right\Vert
_{L^{2}}^{2}  \notag \\
&=&-\int_{\mathbb{R}^{3}}\theta e_{3}\cdot \Delta udx-\sum_{i=1}^{3}\int_{%
\mathbb{R}^{3}}(\partial _{i}u\cdot \nabla )u\partial
_{i}udx-\sum_{i=1}^{3}\int_{\mathbb{R}^{3}}(\partial _{i}u\cdot \nabla
)\theta \partial _{i}\theta dx  \notag \\
&:&=I_{1}+I_{2}+I_{3}.  \label{eq3.1}
\end{eqnarray}%
In the following, we estimate each term at the right-hand side of (\ref%
{eq3.1}) separately below.

To bound $I_{1}$, we integrate by parts and apply H\"{o}lder's inequality to
obtain%
\begin{equation*}
\left\vert I_{1}\right\vert \leq C\left\Vert \nabla u\right\Vert
_{L^{2}}\left\Vert \nabla \theta \right\Vert _{L^{2}}\leq C(\left\Vert
\nabla u\right\Vert _{L^{2}}^{2}+\left\Vert \nabla \theta \right\Vert
_{L^{2}}^{2}).
\end{equation*}

In order to deal with the terms $I_{2}$ and $I_{3}$, we need the following
elegant Machihara-Ozawa inequality \cite{MO} (see also  Meyer \cite{Mey})%
\begin{equation}
\Vert \nabla u\Vert _{L^{4}}^{2}\leq C\left\Vert u\right\Vert _{\dot{B}%
_{\infty ,\infty }^{0}}\left\Vert \Delta u\right\Vert _{L^{2}}.  \label{eq27}
\end{equation}%
We now bound $I_{2}$. By (\ref{eq27}) and Young's inequality
\begin{eqnarray*}
\left\vert I_{2}\right\vert &\leq &C\left\Vert \nabla u\right\Vert
_{L^{2}}\left\Vert \nabla u\right\Vert _{L^{4}}^{2} \\
&\leq &C\left\Vert \nabla u\right\Vert _{L^{2}}\left\Vert u\right\Vert _{%
\dot{B}_{\infty ,\infty }^{0}}\left\Vert \Delta u\right\Vert _{L^{2}} \\
&\leq &C\left\Vert \nabla u\right\Vert _{L^{2}}\left\Vert u\right\Vert
_{BMO}\left\Vert \Delta u\right\Vert _{L^{2}} \\
&\leq &\frac{1}{2}\left\Vert \Delta u\right\Vert _{L^{2}}^{2}+C\left\Vert
\nabla u\right\Vert _{L^{2}}^{2}\left\Vert u\right\Vert _{BMO}^{2}.
\end{eqnarray*}%
By integration by parts, we can rewrite and estimate $I_{3}$ as follows%
\begin{eqnarray*}
\left\vert I_{3}\right\vert &=&\left\vert \sum_{i=1}^{3}\int_{\mathbb{R}%
^{3}}(\partial _{i}u\cdot \nabla )\theta \cdot \partial _{i}\theta
dx\right\vert =\left\vert \sum_{i,j,k=1}^{3}\int_{\mathbb{R}^{3}}\partial
_{i}(\partial _{i}\theta _{k}\partial _{k}\theta _{j})u_{j}dx\right\vert \\
&\leq &C\left\Vert u\right\Vert _{BMO}\left\Vert \nabla \theta \right\Vert
_{L^{2}}\left\Vert \Delta \theta \right\Vert _{L^{2}} \\
&\leq &\frac{1}{6}\left\Vert \Delta \theta \right\Vert
_{L^{2}}^{2}+C\left\Vert u\right\Vert _{BMO}^{2}\left\Vert \nabla \theta
\right\Vert _{L^{2}}^{2}.
\end{eqnarray*}%
Combining the estimates for $I_{1}$, $I_{2}$ and $I_{3}$, we find%
\begin{eqnarray*}
&&\frac{d}{dt}\left( \left\Vert \nabla u\right\Vert _{L^{2}}^{2}+\left\Vert
\nabla \theta \right\Vert _{L^{2}}^{2}\right) +\left\Vert \Delta
u\right\Vert _{L^{2}}^{2}+\left\Vert \Delta \theta \right\Vert _{L^{2}}^{2}
\\
&\leq &C(1+\left\Vert u\right\Vert _{BMO}^{2})(\left\Vert \nabla
u\right\Vert _{L^{2}}^{2}+\left\Vert \nabla \theta \right\Vert _{L^{2}}^{2}).
\end{eqnarray*}%
Using the Gronwall  inequality on the time interval $[T_{0},t]$, one has
the following inequality%
\begin{eqnarray*}
&&\left\Vert \nabla u(\cdot ,t)\right\Vert _{L^{2}}^{2}+\left\Vert \nabla
\theta (\cdot ,t)\right\Vert _{L^{2}}^{2}+\int_{T_{0}}^{t}(\left\Vert \Delta
u(\cdot ,\tau )\right\Vert _{L^{2}}^{2}+\left\Vert \Delta \theta (\cdot
,\tau )\right\Vert _{L^{2}}^{2})d\tau \\
&\leq &\left( \left\Vert \nabla u(\cdot ,T_{0})\right\Vert
_{L^{2}}^{2}+\left\Vert \nabla \theta (\cdot ,T_{0})\right\Vert
_{L^{2}}^{2}\right) \exp \left( C\int_{T_{0}}^{t}\left\Vert u(\cdot ,\tau
)\right\Vert _{BMO}^{2}d\tau \right) .
\end{eqnarray*}%
Let us denote for any $t\in \lbrack T_{0},T)$,
\begin{equation}
F(t)\triangleq \underset{T_{0}\leq \tau \leq t}{\max }\left( \left\Vert
u(\cdot ,\tau )\right\Vert _{H^{2}}^{2}+\left\Vert \theta (\cdot ,\tau
)\right\Vert _{H^{2}}^{2}\right) .  \label{eq52}
\end{equation}%
It should be noted that the function $F(t)$ is nondecreasing. Using (\ref%
{eq21}), we obtain
\begin{eqnarray*}
&&\left\Vert \nabla u(\cdot ,t)\right\Vert _{L^{2}}^{2}+\left\Vert \nabla
\theta (\cdot ,t)\right\Vert _{L^{2}}^{2}+\int_{T_{0}}^{t}(\left\Vert \Delta
u(\cdot ,\tau )\right\Vert _{L^{2}}^{2}+\left\Vert \Delta \theta (\cdot
,\tau )\right\Vert _{L^{2}}^{2})d\tau \\
&\leq &C(T_{0})\exp \left( C\int_{T_{0}}^{t}(1+\left\Vert u(\cdot ,\tau
)\right\Vert _{\dot{B}_{\infty ,\infty }^{0}}^{2}\log (\left\Vert u(\cdot
,\tau )\right\Vert _{H^{2}}+\left\Vert \theta (\cdot ,\tau )\right\Vert
_{H^{2}}))d\tau \right) \\
&\leq &C(T_{0})\exp \left( C\int_{T_{0}}^{t}\left\Vert \nabla u(\cdot ,\tau
)\right\Vert _{\dot{B}_{\infty ,\infty }^{-1}}^{2}\log (\left\Vert u(\cdot
,\tau )\right\Vert _{H^{2}}^{2}+\left\Vert \theta (\cdot ,\tau )\right\Vert
_{H^{2}}^{2})d\tau \right) \\
&\leq &C(T_{0})\exp \left( C\int_{T_{0}}^{t}\left\Vert \nabla u(\cdot ,\tau
)\right\Vert _{\overset{.}{B}_{\infty ,\infty }^{-1}}^{2}d\tau \underset{%
T_{0}\leq \tau \leq t}{\sup }\log (\left\Vert u(\cdot ,\tau )\right\Vert
_{H^{2}}^{2}+\left\Vert \theta (\cdot ,\tau )\right\Vert _{H^{2}}^{2})\right)
\\
&\leq &C(T_{0})\exp \left( C\int_{T_{0}}^{t}\left\Vert \nabla u(\cdot ,\tau
)\right\Vert _{\overset{.}{B}_{\infty ,\infty }^{-1}}^{2}d\tau \log \underset%
{T_{0}<\tau \leq t}{\sup }(\left\Vert u(\cdot ,\tau )\right\Vert
_{H^{2}}^{2}+\left\Vert \theta (\cdot ,\tau )\right\Vert _{H^{2}}^{2})\right)
\\
&\leq &C(T_{0})\exp \left( C\epsilon \log F(t)\right) \\
&\leq &C(T_{0})\left[ F(t)\right] ^{C\epsilon },
\end{eqnarray*}%
where
\begin{equation*}
C(T_{0})=C\left( \left\Vert \nabla u(\cdot ,T_{0})\right\Vert
_{L^{2}}^{2}+\left\Vert \nabla \theta (\cdot ,T_{0})\right\Vert
_{L^{2}}^{2}\right) .
\end{equation*}

Next, applying $\Delta $ to the equations (\ref{eq1.1})$_{1}$, (\ref{eq1.1})$%
_{2}$, taking the $L^{2}$ inner product of the obtained equations with $%
-\Delta u$ and $-\Delta \theta $, respectively, adding them up and using the
incompressible conditions $\nabla \cdot u=0$, we arrive at%
\begin{eqnarray}
&&\frac{1}{2}\frac{d}{dt}(\left\Vert \Delta u\right\Vert
_{L^{2}}^{2}+\left\Vert \Delta \theta \right\Vert _{L^{2}}^{2})+\left\Vert
\nabla ^{3}u\right\Vert _{L^{2}}^{2}+\left\Vert \nabla ^{3}\theta
\right\Vert _{L^{2}}^{2}  \notag \\
&=&\int_{\mathbb{R}^{3}}\Delta (\theta e_{3})\cdot \Delta udx-\int_{\mathbb{R%
}^{3}}\Delta (u\cdot \nabla u)\cdot \Delta udx-\int_{\mathbb{R}^{3}}\Delta
(u\cdot \nabla \theta )\cdot \Delta \theta dx  \notag \\
&\leq &\left\vert \int_{\mathbb{R}^{3}}\Delta (\theta e_{3})\cdot \Delta
udx\right\vert +\left\vert \int_{\mathbb{R}^{3}}(\Delta u\cdot \nabla
u)\cdot \Delta udx\right\vert +2\sum_{i=1}^{3}\left\vert \int_{\mathbb{R}%
^{3}}(\partial _{i}u\cdot \nabla \partial _{i}u)\cdot \Delta udx\right\vert
\notag \\
&&+\left\vert \int_{\mathbb{R}^{3}}(\Delta u\cdot \nabla \theta )\cdot
\Delta \theta dx\right\vert +2\sum_{i=1}^{3}\left\vert \int_{\mathbb{R}%
^{3}}(\partial _{i}u\cdot \nabla \partial _{i}\theta )\cdot \Delta \theta
dx\right\vert  \notag \\
&=&\sum\limits_{k=1}^{5}A_{k}.  \label{eq44}
\end{eqnarray}%
Now we will estimate the terms on the right-hand side of (\ref{eq44}) one by
one as follows. Let us begin with estimating the term $A_{1}$.

Using Lemma \ref{lem3} with $p=q=r=j=2$, $k=3$ and $m=1$, $A_{1}$ can be
bounded above as follows:%
\begin{eqnarray*}
A_{1} &\leq &C\left\Vert \Delta u\right\Vert _{L^{2}}\left\Vert \Delta
\theta \right\Vert _{L^{2}} \\
&\leq &C\left\Vert \nabla u\right\Vert _{L^{2}}^{\frac{1}{2}}\left\Vert
\nabla ^{3}u\right\Vert _{L^{2}}^{\frac{1}{2}}\left\Vert \nabla \theta
\right\Vert _{L^{2}}^{\frac{1}{2}}\left\Vert \nabla ^{3}\theta \right\Vert
_{L^{2}}^{\frac{1}{2}} \\
&=&\left( \left\Vert \nabla ^{3}u\right\Vert _{L^{2}}^{2}\right) ^{\frac{1}{4%
}}\left( \left\Vert \nabla ^{3}\theta \right\Vert _{L^{2}}^{2}\right) ^{%
\frac{1}{4}}\left( C\left\Vert \nabla u\right\Vert _{L^{2}}\left\Vert \nabla
\theta \right\Vert _{L^{2}}\right) ^{\frac{1}{2}} \\
&\leq &\frac{1}{16}\left\Vert \nabla ^{3}u\right\Vert _{L^{2}}^{2}+\frac{1}{%
16}\left\Vert \nabla ^{3}\theta \right\Vert _{L^{2}}^{2}+C\left\Vert \nabla
u\right\Vert _{L^{2}}\left\Vert \nabla \theta \right\Vert _{L^{2}} \\
&\leq &\frac{1}{16}\left\Vert \nabla ^{3}u\right\Vert _{L^{2}}^{2}+\frac{1}{%
16}\left\Vert \nabla ^{3}\theta \right\Vert _{L^{2}}^{2}+C\left( \left\Vert
\nabla u\right\Vert _{L^{2}}^{2}+\left\Vert \nabla \theta \right\Vert
_{L^{2}}^{2}\right) .
\end{eqnarray*}%
Let us now recall Gagliardo--Nirenberg's inequality%
\begin{equation*}
\left\Vert \Delta f\right\Vert _{L^{4}}\leq C\left\Vert \nabla f\right\Vert
_{L^{2}}^{\frac{1}{8}}\left\Vert \nabla ^{3}f\right\Vert _{L^{2}}^{\frac{7}{8%
}}.
\end{equation*}%
Thus, we obtain
\begin{eqnarray*}
A_{2},A_{3} &\leq &C\left\Vert \nabla u\right\Vert _{L^{2}}\left\Vert \Delta
u\right\Vert _{L^{4}}^{2} \\
&\leq &C\left\Vert \nabla u\right\Vert _{L^{2}}\left\Vert \nabla
u\right\Vert _{L^{2}}^{\frac{1}{4}}\left\Vert \nabla ^{3}u\right\Vert
_{L^{2}}^{\frac{7}{4}} \\
&=&C\left\Vert \nabla u\right\Vert _{L^{2}}^{\frac{5}{4}}\left\Vert \nabla
^{3}u\right\Vert _{L^{2}}^{\frac{7}{4}}=\left( C\left\Vert \nabla
u\right\Vert _{L^{2}}^{10}\right) ^{\frac{1}{8}}\left( \left\Vert \nabla
^{3}u\right\Vert _{L^{2}}^{2}\right) ^{\frac{7}{8}} \\
&\leq &\frac{1}{16}\left\Vert \nabla ^{3}u\right\Vert
_{L^{2}}^{2}+C\left\Vert \nabla u\right\Vert _{L^{2}}^{10}.
\end{eqnarray*}%
Similarly to the estimate of $A_{1}$, the terms $A_{4}$ and $A_{5}$ can be
bounded  above as
\begin{eqnarray*}
A_{4},A_{5} &\leq &C\left\Vert \nabla \theta \right\Vert _{L^{2}}\left\Vert
\Delta \theta \right\Vert _{L^{4}}\left\Vert \Delta u\right\Vert _{L^{4}} \\
&\leq &C\left\Vert \nabla \theta \right\Vert _{L^{2}}\left( \left\Vert
\Delta u\right\Vert _{L^{4}}^{2}+\left\Vert \Delta \theta \right\Vert
_{L^{4}}^{2}\right) \\
&\leq &C\left\Vert \nabla \theta \right\Vert _{L^{2}}\left\Vert \nabla
u\right\Vert _{L^{2}}^{\frac{1}{4}}\left\Vert \nabla ^{3}u\right\Vert
_{L^{2}}^{\frac{7}{4}}+C\left\Vert \nabla \theta \right\Vert _{L^{2}}^{\frac{%
5}{4}}\left\Vert \nabla ^{3}\theta \right\Vert _{L^{2}}^{\frac{7}{4}} \\
&\leq &\frac{1}{4}\left\Vert \nabla ^{3}u\right\Vert
_{L^{2}}^{2}+C\left\Vert \nabla \theta \right\Vert _{L^{2}}^{8}\left\Vert
\nabla u\right\Vert _{L^{2}}^{2}+\frac{1}{2}\left\Vert \nabla ^{3}\theta
\right\Vert _{L^{2}}^{2}+C\left\Vert \nabla \theta \right\Vert _{L^{2}}^{10}
\\
&\leq &\frac{1}{16}\left\Vert \nabla ^{3}u\right\Vert _{L^{2}}^{2}+\frac{1}{4%
}\left\Vert \nabla ^{3}\theta \right\Vert _{L^{2}}^{2}+C\left\Vert \nabla
\theta \right\Vert _{L^{2}}^{8}\left( \left\Vert \nabla u\right\Vert
_{L^{2}}^{2}+\left\Vert \nabla \theta \right\Vert _{L^{2}}^{2}\right) .
\end{eqnarray*}%
Summarizing all the estimates and absorbing the dissipative term, we can
derive
\begin{eqnarray*}
\frac{d}{dt}(\left\Vert \Delta u\right\Vert _{L^{2}}^{2}+\left\Vert \Delta
\theta \right\Vert _{L^{2}}^{2}) &\leq &C\left\Vert \nabla u\right\Vert
_{L^{2}}^{10}+C\left\Vert \nabla \theta \right\Vert _{L^{2}}^{8}\left(
\left\Vert \nabla u\right\Vert _{L^{2}}^{2}+\left\Vert \nabla \theta
\right\Vert _{L^{2}}^{2}\right) \\
&\leq &C\left( \left\Vert \nabla u\right\Vert _{L^{2}}^{8}+\left\Vert \nabla
\theta \right\Vert _{L^{2}}^{8}\right) \left( \left\Vert \nabla u\right\Vert
_{L^{2}}^{2}+\left\Vert \nabla \theta \right\Vert _{L^{2}}^{2}\right) \\
&\leq &C\left( \left\Vert \nabla u\right\Vert _{L^{2}}^{2}+\left\Vert \nabla
\theta \right\Vert _{L^{2}}^{2}\right) ^{4}\left( \left\Vert \nabla
u\right\Vert _{L^{2}}^{2}+\left\Vert \nabla \theta \right\Vert
_{L^{2}}^{2}\right) \\
&\leq &C\left( \left\Vert \nabla u\right\Vert _{L^{2}}^{2}+\left\Vert \nabla
\theta \right\Vert _{L^{2}}^{2}\right) ^{5} \\
&\leq &C(T_{0})\left[ F(t)\right] ^{5C\epsilon }.
\end{eqnarray*}%
Integrating the above estimate over interval $(T_{0},t)$ and observing that $%
F(t)$ is a monotonically increasing function, we thus have
\begin{eqnarray*}
&&\left\Vert \Delta u(\cdot ,t)\right\Vert _{L^{2}}^{2}+\left\Vert \Delta
\theta (\cdot ,t)\right\Vert _{L^{2}}^{2} \\
&\leq &\left\Vert \Delta u(\cdot ,T_{0})\right\Vert _{L^{2}}^{2}+\left\Vert
\Delta \theta (\cdot ,T_{0})\right\Vert
_{L^{2}}^{2}+C(T_{0})\int_{T_{0}}^{t} \left[ F(\tau )\right] ^{5C\epsilon
}d\tau .
\end{eqnarray*}%
By using (\ref{eq52}), it follows that%
\begin{eqnarray*}
F(t) &\leq &\left\Vert u(\cdot ,T_{0})\right\Vert _{H^{2}}^{2}+\left\Vert
\theta (\cdot ,T_{0})\right\Vert _{H^{2}}^{2}+C\int_{T_{0}}^{t}\left[ F(\tau
)\right] ^{5C\epsilon }d\tau \\
&\leq &\left\Vert u(\cdot ,T_{0})\right\Vert _{H^{2}}^{2}+\left\Vert \theta
(\cdot ,T_{0})\right\Vert _{H^{2}}^{2}+C(T_{0})(t-T_{0})\left[ F(t)\right]
^{5C\epsilon }.
\end{eqnarray*}%
Choosing $\epsilon $ such that $5C\epsilon <1$, the above
inequality yields for any $t\in \lbrack T_{0},T)$
\begin{equation*}
F(t)\leq C<\infty ,
\end{equation*}%
which implies that $(u,\theta )\in L^{\infty }(0,T;H^{1}(\mathbb{R}%
^{3}))\cap L^{2}(0,T;H^{2}(\mathbb{R}^{3}))$. This completes the proof of
Theorem \ref{th1}.

\begin{remark}
Comparing our result with \cite{Ye}, we have  simplified the proof of Theorem 1.1
in \cite{Ye}, in fact we only need $H^{2}$ a priori estimates of solutions.
\end{remark}

\section{Acknowledgment.}

Part of the work was carried out while the third author was long-term
visitor at University of Catania. The hospitality of Catania University is
graciously acknowledged. This research is partially supported by Piano della
Ricerca 2016-2018 - Linea di intervento 2: "Metodi variazionali ed equazioni
differenziali". M.A. Ragusa wish to thank the support of "RUDN University
Program 5-100".
The authors wish to express their thanks to the referees for
their very careful reading of the paper, giving valuable comments and
helpful suggestions.


\end{document}